\documentclass[11pt]{article}
\usepackage{amsmath,amsthm}
\usepackage{pstricks,pst-node,pst-tree}
\usepackage{latexsym}
\usepackage{amssymb,amscd}
\usepackage{geometry, graphicx}
\usepackage{url}
\numberwithin{equation}{section}
\let \:=\colon
\let \beg=\begin

\let \hra=\hookrightarrow
\let \mb=\mathbb
\let \mc= \mathcal

\let \ra=\rightarrow

\let \ka=\kappa
\let \la=\lambda
\let \La=\Lambda

\let \al=\alpha

\let \de=\delta

\let \ep=\epsilon

\let \fl=\flushleft
\let \fr=\frac
\let \mc=\mathcal

\let \part=\partial

\let \sig=\sigma
\let \sub=\subset
\let \wh=\widehat

\beg{document}
\title{\bf Rational curves of degree 16 on a general heptic fourfold}
\author{\small \hspace{-40pt} Ethan Cotterill \\ 
\small Departamento de Matem\'atica, Universidade de Coimbra, Apartado 3008, 3001-454 Coimbra, Portugal.
}
\date{\empty}
\maketitle

\beg{abstract}
According to a conjecture of H. Clemens, the dimension of the space of rational curves on a general projective hypersurface should equal the number predicted by a na\"ive dimension count. In the case of a general hypersurface of degree 7 in $\mb{P}^5$, the conjecture predicts that the only rational curves should be lines. This has been verified by Hana and Johnsen for rational curves of degree at most 15. Here we extend their results to show that no rational curves of degree 16 lie on a general heptic fourfold.
\end{abstract}

\section{Introduction}
A celebrated conjecture of H. Clemens predicts that the number of smooth rational curves of fixed degree on a general quintic threefold in $\mb{P}^4$ is finite. General here  means outside of a countable union of Zariski-closed proper subvarieties of the projective space of hypersurfaces, i.e., what is sometimes referred to as ``very general". In \cite{C1} and \cite{C2}, we proved a stronger statement for curves of degree at most 11; namely, that the corresponding incidence scheme of curves in quintics is irreducible, and that, the only singular rational curves of degree at most 11 on a general quintic are six-nodal plane quintics. In both cases our arguments were based on involved analyses of {\it initial ideals} associated to rational curves in $\mb{P}^4$, and their cohomology.

{\fl In} this note, we will apply the same technique to the study of rational curves on a general hypersurface of degree 7 in $\mb{P}^5$. A na\"ive dimension count leads one to expect that the only rational curves on such a fourfold should be lines. Indeed, Hana and Johnsen \cite{HJ} used estimates of ideal sheaf cohomology based on the method of \cite{C1,C2} to verify that the only rational curves of degree at most 15 are lines.

{\fl We} obtain the following extension of the result of \cite{HJ}:

\beg{thm}\label{mainthm}
The only rational curves of degree at most 16 on a general heptic fourfold are lines.
\end{thm}

{\fl We} prove Theorem~\ref{mainthm} first for {\it nondegenerate} rational curves, i.e. those not contained in hyperplanes.

\section{Nondegenerate curves and restricted tangent bundles}
Every irreducible, nondegenerate degree-$d$ rational curve in $\mb{P}^n$ is the image of a degree-$d$ morphism $f:\mb{P}^1 \ra \mb{P}^n$; as such, it determines a point $[f]$ of the variety $M^d_n:= \mbox{Mor}_d(\mb{P}^1,\mb{P}^n)$. Concretely, $M^d_n$ is a space of $(n+1)$-tuples of polynomials in 2 homogeneous parameters for $\mb{P}^1$, modulo the action of $\mbox{PGL}(2)$ and scalar multiplication. Thus $M_n^d$ is a smooth scheme of dimension $(n+1)d-4$.

{\fl The} tangent space of $M^d_n$ at a point $[f]$ may be canonically identified with the space of sections $H^0(\mb{P}^1,f^* \mc{T}_{\mb{P}^n})$ of the restricted tangent bundle. Verdier \cite{V} and Ramella \cite{R} showed that the stratum of the mapping space associated to a particular splitting type of $f^* \mc{T}_{\mb{P}^n}$ is smooth and of the expected codimension. It is convenient to index the strata by the $n$-tuples $(a_1,\dots,a_n)$ associated to the twisted decompositions
\[
f^* (\mc{T}_{\mb{P}^n}(1))= \mc{O}_{\mb{P}^1}(a_1) \oplus \cdots \oplus \mc{O}_{\mb{P}^1}(a_n)
\]
such that $\sum_{i=1}^n a_i=d$. The $(a_1,\dots,a_n)$-stratum is then of codimension equal to \\ $\sum_{i \neq j; a_i \geq a_j} \max\{0,a_i-a_j-1\}$.

\subsection{Nondegenerate rational curves on heptics in $\mb{P}^5$}\label{2.1}
Let $\mc{I}_{d,g,i}$ denote the incidence variety of irreducible, nondegenerate rational curves $C$ of degree $d$ and arithmetic genus $g$ on heptics for which $h^1(\mb{P}^5, \mc{I}_{C/\mb{P}^5}(7))=i$). Let $\mb{P}^N$ denote the projective space of heptic fourfolds. A standard argument (see, e.g., \cite{C1}) implies that the refined incidence variety is of dimension
\[
\mbox{dim}(\mc{I}_{d,g,i})= N+1-d+g+i.
\]
{\fl In} particular, $\mc{I}_{d,g,i}$ fails to dominate the space of heptics whenever $g+i \leq d-2$.

{\fl A} further useful refinement of the incidence variety is via the stratification by splitting type of the restricted tangent bundle. Given a nondegenerate degree-$d$ morphism $f: \mb{P}^1 \ra \mb{P}^5$, the result \cite[Prop 1.2]{GLP} implies that the image $C$ of $f$ is $8$-regular whenever the corresponding splitting type $(a_1,a_2,a_3,a_4,a_5)$ satisfies $a_i+a_j \leq 8$ for all $i \neq j \in \{1,2,3,4,5\}$.

\subsection{A stratification scheme for cohomological estimates}\label{strat_sch}
The result of \cite{GLP} cited in the preceded subsection implies that the only morphisms $f$ of degree 16 with $h^1(\mb{P}^5, \mc{I}_{C/\mb{P}^5}(7)) \neq 0$ belong to strata $(a_1,a_2,a_3,a_4,a_5)$ with
\beg{equation}\label{basic_restriction}
a_1 \geq a_2 \geq a_3 \geq a_4 \geq a_5, \text{ and } a_1+ a_2 \geq 9.
\end{equation}
{\fl Note} that the most generic such stratum is $(5,4,3,2,2)$, and is of codimension 7 in the mapping space.

{\fl To} show that the image $C$ of such a morphism $f$ fails to lie on a general heptic, it suffices to show that
\beg{equation}\label{cohest}
g+i \leq 14+ \sum_{i \neq j; a_i \geq a_j} \max\{0,a_i-a_j-1\}.
\end{equation}
where $i:=h^1(\mb{P}^5, \mc{I}_{C/\mb{P}^5}(7))$.

\subsection{Bounds on $g+i$ via generic initial ideals and Castelnuovo's estimate}\label{gins}

Given any homogeneous ideal $I$ in $(n+1)$ variables $x_i, 0 \leq i \leq n$,
together with a
choice of
partial order $<$ on the monomials $m(x_i)$, the {\it initial ideal} of $I$ with respect to $<$ 
is generated by leading terms of elements in $I$. Given any subscheme $C \sub \mathbb{P}^n$, upper-semicontinuity implies that the regularity and values $h^j$
of the ideal sheaf $\mathcal{I}_C$ are majorized by those of the initial ideal sheaf
$\mbox{in}(\mathcal{I}_C)$ with respect to any partial order on the
monomials of $\mathbb{P}^n$. Moreover, Macaulay's theorem establishes the Hilbert
functions of $\mathcal{I}_C$ and $\mbox{in}(\mathcal{I}_C)$ agree. Replacing $C$ by
a {\it general} $\mbox{PGL}(n+1)$-translate, we obtain {\it generic initial
  ideal}
of $C$, written $\mbox{gin}(\mathcal{I}_C$). 

{\fl There} is a unique generic initial ideal associated to any
  ideal. Moreover, the following three statements are equivalent (see \cite[Thm.~2.30]{Gr}).
\begin{enumerate}
\item The monomial ideal $\mbox{gin }\mathcal{I}_C$ is saturated, in
  the ideal-theoretic sense. In other words, for all $m \in (x_0,\dots,x_n)$,
\vspace{-.2cm}
\begin{equation*}
g \cdot m \in  \mbox{gin
  }\mathcal{I}_C \Rightarrow g \in \mbox{gin
  }\mathcal{I}_C.
\vspace{-.25cm}
\end{equation*}
\item The homogeneous ideal $I_C$
  is saturated.
\item No minimal generator of $\mbox{gin
  }\mathcal{I}_C \subset \mathbb{C}[x_0,\dots,x_n]$ is divisible by $x_n$.
\end{enumerate}

{\fl Now} let $I \subset \mathbb{C}[x_0,\dots,x_n]$ be any ideal that
is {\it Borel-fixed}, i.e., fixed under the action of upper triangular
matrices $\mathcal{T} \subset \mbox{PGL}(n+1)$; generic initial ideals are always of this type, and $\mathcal{T}$-fixedness
means exactly that $I$ is minimally generated by monomials $P$ for which
\[
x_i/x_j \cdot P \in I \text{ whenever } i<j.
\]

{\fl Then} $I$ is
saturated if and only if none of its minimal generators is divisible by $x_n$.
Hereafter, we work exclusively
with {\it saturated} generic initial ideals. An important fact is that {\it the Castelnuovo-Mumford regularity of any saturated ideal sheaf $\mc{I}$ equals that of its associated generic initial ideal sheaf, which in turn equals the maximum degree of any minimal generator of $\mbox{gin}(I)$.}

{\fl The} basic two-step procedure for bounding $g+i$ developed in my earlier paper \cite[Sect. $1$]{C1} is as follows.

\beg{itemize}
\item[(1)] Determine possible monomial ideals arising as gins associated to general hyperplane sections of rational curves $C$.
\item[(2)] Determine possible curve gins corresponding to each hyperplane section gin, and compute the corresponding values of $g+i$.
\end{itemize}

{\fl Indeed}, part (1) plays a prominent r\^ole in \cite{C1} and \cite{C2}. However, the bounds on cohomology required for Clemens' conjecture for rigidity of rational curves on a general quintic threefold are more restrictive than \eqref{cohest} (here, we have $d$ more degrees of freedom!). Instead of exhaustively listing Borel-fixed monomial ideals, the following observations will suffice.

\beg{itemize}
\item A celebrated result of Castelnuovo establishes that the arithmetic genus $g$ of any irreducible, nondegenerate curve $C \sub \mb{P}^n$ satisfies
\beg{equation}\label{castelnuovo_bound}
g(C) \leq \binom{m}{2}(n-1)+m \ep
\end{equation}
where $m= \lfloor \fr{d-1}{n-1} \rfloor$ and $d-1=m(n-1)+\ep$. In particular, when $d=16$ and $n=5$, we find $g(C) \leq 21$. An important remark is that {\it the proof of Castelnuovo's result uses only the fact that a general hyperplane section of $C$ is a set of points in uniform position}; see, e.g., \cite{Ha}.
\item The generic initial ideal, or {\it gin}, of an irreducible, nondegenerate rational curve $C \sub \mb{P}^5$ is related to the gin of its generic hyperplane section $\La$ by a series of {\it rewriting rules} that (immediately) generalize those given in \cite[Sect. 1]{C1}; see \cite[Lem. 3.6]{HJ}, and the arithmetic genus $g$ and cohomology $i$ of $C$ verify
\beg{equation}\label{basicineq}
g+i \leq g_{\La}
\end{equation}
where $g_{\La}$ is the arithmetic genus of the cone $C_{\La}$ with vertex $(0,0,0,0,0,1)$ over $\La \sub \mb{P}^4$, with respect to the inclusion of $\mb{P}^4$ in $\mb{P}^5$ given by the inclusion of sets of homogeneous variables $\{x_i\}_{i=0,\dots,4} \sub \{x_i\}_{i=0,\dots,5}$.
\end{itemize}

\subsection{Proof of Theorem~\ref{mainthm} for nondegenerate curves}
Let $C$ denote an irreducible, nondegenerate rational curve of degree 16 in $\mb{P}^5$. Let $\La$ denote a generic hyperplane section of $C$, and let $C_{\La}$ denote the cone over $\La$, as defined in the previous subsection. While the curve $C_{\La}$ is itself reducible, its generic hyperplane is a set of points in uniform position, so it also satisfies Castelnuovo's genus bound! Consequently, combining \eqref{castelnuovo_bound} with \eqref{basicineq}, we deduce
\beg{equation}\label{crude_estimate}
g+i \leq 21.
\end{equation}
Now compare \eqref{crude_estimate} with the required estimate \eqref{cohest}. In Subsection~\ref{strat_sch}, we noted that if $i \neq 0$, then in fact $C$ belongs to a splitting stratum of codimension (at least) 7, so in this case \eqref{cohest} follows immediately and we conclude! 

{\fl Without} loss of generality, therefore, we may assume $i=0$, and $15 \leq g \leq 21$. The corresponding rational curves are highly singular, and we defer handling them until the final section of the paper.

\section{Degenerate rational curves}
Here we adopt a strategy analogous to the one used in the preceding section to prove Theorem~\ref{mainthm} for curves that linearly span $\mb{P}^4$ or $\mb{P}^3$, respectively (B\'ezout's theorem automatically precludes plane curves of degree 16 from lying on a heptic).

\subsection{Proof of Theorem~\ref{mainthm} for curves that span $\mb{P}^4$}
We begin by assembling the required ingredients, namely the analogue of \eqref{cohest}, and the crude estimate \eqref{crude_estimate} coming from Castelnuovo's genus bound.

{\fl This time}, the required estimate, which replaces \eqref{cohest}, is
\beg{equation}\label{cohest2}
g+i \leq 28 + \sum_{i \neq j; a_i \geq a_j} \max\{0,a_i-a_j-1\}
\end{equation}
where $(a_i)_{1 \leq i \leq 4}$ gives the splitting type of the (twisted) restricted tangent bundle arising from the $\mb{P}^4$ in which our curve $C$ is embedded. For the summand 28, see, e.g. the proof of \cite[Lem 2.4]{HJ}. On the other hand,
\beg{equation}\label{crude_estimate2}
g+i \leq 30
\end{equation}
is the Castelnuovo estimate.

{\fl For} the sake of argument, assume $i \neq 0$. Then \cite[Prop 1.2]{GLP} implies that our curve $C$ belongs to the closure of the codimension-one splitting stratum $(5,4,4,3)$. However, the splitting stratum is irreducible, and a explicit verification carried out in Macaulay2 \cite{GS} in the last section of the paper shows that a generic curve belonging to the stratum has $i=0$. Consequently, whenever $i \neq 0$, the required estimate \eqref{cohest2} may be weakened so that it becomes precisely the estimate arising from Castelnuovo's genus bound, and we conclude.

{\fl Without} loss of generality, therefore, we may assume $i=0$, and $28 \leq g \leq 30$. The corresponding rational curves are (highly) singular, and to conclude it suffices to note that singular nondegenerate rational curves in $\mb{P}^4$ determine a codimension-2 sublocus of $M^d_4$ for all $d \geq 5$.

\subsection{Proof of Theorem~\ref{mainthm} for curves that span $\mb{P}^3$}
In this case, the estimate we need is
\beg{equation}\label{cohest3}
g+i \leq 42,
\end{equation}
while Castelnuovo's genus bound yields
\beg{equation}\label{crude_estimate2}
g+i \leq 49.
\end{equation}

{\fl In this case,} we are forced to consider possible gins of general hyperplane sections $\La$ of the corresponding curves. As explained in \cite[Sect. 1]{C1}, we have
\beg{equation}\label{mingens}
\mbox{gin}(\mc{I}_{\La})=(x_0^k,x_0^{k-1}x_1^{\la_{k-1}}, \dots, x_0x_1^{\la_1},x_1^{\la_0})
\end{equation}
for suitably chosen positive integers $\la_i, 0 \leq i \leq k-2$ which moreover satisfy
\begin{equation}\label{GP}
\la_i-2 \leq \la_{i+1} \leq \la_{i-1}
\end{equation}
by a result of Gruson and Peskine \cite{GP}. In \cite{C1} we describe a natural scheme for representing a set of minimal generators for $\mbox{gin}(\mc{I}_{\La})$ as a tree generated by {\it rewriting rules}; the number of rewritings is precisely
\[
\deg(\La)= \sum_{i=0}^{k-2} \la_i. 
\]
Moreover, a result of Ballico \cite{Ba} for schemes of points in uniform position implies $\mbox{gin}(\mc{I}_{\La})$ is 9-regular. The combinatorial genus formula of \cite[Lem 1.5.2]{C1} ($\binom{m+4}{4}$ needs to be replaced by $\binom{m+3}{3}$, and $10m$ by $16m$, as $C$ is of degree 16 and spans $\mb{P}^3$) now yields
\beg{equation}\label{glambdafmla}
\beg{split}
g_{\La}&= 16m+1- \binom{m+3}{3}+ h^0(\mc{I}_{C_{\La}}(m)) \text{ with } m=9\\
&= 16m+1- \binom{m+3}{3}+ \binom{m+3-k}{3}+ \sum_{i=0}^{k-1}\binom{m+2-(i+\la_{i})}{2}.
\end{split}
\end{equation}

{\fl Combining} \eqref{GP} and \eqref{glambdafmla}, we deduce:
\beg{itemize}
\item $g_{\La} \leq 31$ whenever $k \geq 3$.
\item If $k=2$, then $\mbox{gin}(\mc{I}_{\La})= (x_0^2,x_0x_1^7,x_1^9)$, in which case $g_{\La}=49$.
\end{itemize}
{\fl Therefore, in light} of the required estimate \eqref{cohest3}, we may assume $\mbox{gin}(\mc{I}_{\La})= (x_0^2,x_0x_1^7,x_1^9)$ without loss of generality. However, as explained in \cite[Sect. 2]{C1}, the (tree of minimal generators associated to the) gin of an irreducible curve that linearly spans $\mb{P}^3$ (viewed as a subscheme of the $\mb{P}^3$ it linearly spans) is obtained from the (tree of minimal generators of the) gin of its generic hyperplane section by applying a finite number of {\it $C$-rewriting rules}, as in Table~\ref{table_one}. Moreover, the essential quantity $g+i$ decreases by 1 every time a $C$-rewriting is applied to a monomial of degree less than 7. 

{\fl On} the other hand, rational curves of degree 16 in $\mb{P}^3$ that lie on a quadric determine a sublocus of $M^{16}_3$ of codimension greater than 7, which is the basic disparity between the Castelnuovo estimate and that which is required. Consequently, letting $C$ denote any rational curve with hyperplane gin $(x_0^2,x_0x_1^7,x_1^9)$, we may assume the quadratic generator $x_0^2$ is rewritten, and thereby replace the estimate \eqref{crude_estimate2} by $g+i \leq 48$. On the other hand, B\'ezout's theorem implies that $C$ lies on at most one linearly independent cubic surface in $\mb{P}^3$, and for the same reason it cannot lie on a linearly independent pair of surfaces of degrees $(3,4)$ or even $(4,4)$ (for the latter case, use Max Noether's $AF+BG$ theorem). These considerations force a minimum of {\it nine} $C$-rewritings of monomials of degrees less than 7, which improves our cohomological estimate to $g+i \leq 40$. We conclude that no rational curves of degree 16 that linearly span $\mb{P}^3$ lie on a general heptic.

\begin{table}
\caption{$C$-rules for nondegenerate irreducible curves}\label{table_one}
\vspace{-.3cm}
$$\begin{array}{ll}
\hline
\vspace{0.1cm}
1.&x_0^e
\mapsto x_0^e \cdot(x_0,x_1,x_2), \text{ and }\\
\vspace{.1cm}
2.&x_0^ex_1^f \mapsto x_0^ex_1^f \cdot(x_1,x_2)
\end{array}$$
\end{table}

\section{Loose ends}
\subsection{An explicit calculation in Macaulay2}
In this subsection, we complete the proof that no image of a morphism belonging to $M^{16}_4$ lies on a general heptic, by showing that a general element of the $(5,4,3,3)$ restricted tangent bundle stratum satisfies $i=0$. We follow the same procedure employed in \cite[Pf of Claims 5.1-5.3]{C2}, and use the pullback of the Euler sequence on $\mb{P}^4$ to reverse-engineer an element of the $(5,4,4,3)$ stratum. More precisely, for a generic choice of syzygies of degrees $5,4,4,3$ of the form
\beg{equation}\label{syzygies}
\sum_{0 \leq j \leq 16} g_{i,j} f_{i} =0, 0 \leq i \leq 4
\end{equation}
with $g_{1,j} \in H^0(\mb{P}^1,\mc{O}(5))$; $g_{2,j}, g_{3,j} \in H^0(\mb{P}^1,\mc{O}(4))$; $g_{4,j} \in H^0(\mb{P}^1, \mc{O}(3))$ and
\beg{equation}\label{undetermined_coefs}
f_{i}= \sum_{j=0}^{16} a_{i,j} t^{16-j} u^j \in H^0(\mb{P}^1, \mc{O}(16)),
\end{equation}
the equations \eqref{syzygies} determine a total of $h^0(\mc{O}(21))+ h^0(\mc{O}(20))+ h^0(\mc{O}(20))+ h^0(\mc{O}(19))= 84$ linearly independent relations among the 85 coefficients $a_{i,j}$ coming from \eqref{undetermined_coefs}. In other words, a generic choice of syzygies uniquely determines an element of the splitting stratum. Choosing our syzygies explicitly to be
\beg{equation}\label{explicit_syzygies}
\beg{split}
&t^5 f_0+ t^4u f_1+ (t^3u^2- t^2u^3)f_2+ tu^4 f_3+ u^5 f_4=0, \\
&u^4 f_0+ u^3t f_1+ u^2t^2 f_2+ ut^3 f_3+ t^4 f_4=0, \\
&u^3t f_0+ (u^4-t^4) f_1+ ut^3 f_2+ t^4 f_3+ u^4 f_4=0, \\
&t^3 f_0+ t^2u f_1 + tu^2 f_2+ (t^3+ t^2u+ tu^2+ t^3) f_3+ (t^3-u^3) f_4=0
\end{split}
\end{equation}
we obtain the solution
\[
{\small
\beg{split}
f_0 &=  2024965316t^{15}u -2939521792t^{14}u^2+ 5775517930t^{13}u^3 -3553226765t^{12}u^4 -1878022821t^{11}u^5 \\
&-4727434224t^{10}u^6+ 421222439t^9u^7 -6942771445t^8u^8+ 167619124t^7u^9 -413542814t^6u^{10}\\
&-3988180306t^5u^{11}+ 3858123180t^4u^{12} -1625065825t^3u^{13}+ 1265968623t^2u^{14},\\
f_1 &= -2024965316t^{16}+ 1469760896t^{15}u -1279196964t^{14}u^2+ 657977093t^{13}u^3+ 3046492064t^{12}u^4 -1071413358t^{11}u^5\\
&+ 3131314971t^{10}u^6+ 657977093t^9u^7 -161442087t^8u^8 -2480661437t^7u^9+ 1094059027t^6u^{10}\\
&-4232900355t^5u^{11}+ 2462097431t^4u^{12}+ 607991530t^3u^{13}+ 906871421t^2u^{14}+ 1265968623tu^{15},\\
f_2 &= 1469760896t^{16} -3026560070t^{15}u+ 1893654918t^{14}u^2 -3324744957t^{13}u^3 -1298225797t^{12}u^4 -2145324486t^{11}u^5 \\
&-1485998803t^{10}u^6+ 963954130t^9u^7+ 
1350408778t^8u^8 -1398695708t^7u^9+ 8068973833t^6u^{10}\\
&+ 1146259438t^5u^{11}+  6291363246t^4u^{12}+ 1465865564t^3u^{13} -2841048885t^2u^{14} -547774219tu^{15} -2531937246u^{16} ,\\
f_3 &= -2024965316t^{16}+ 1747363106t^{14}u^2 -1235677825t^{13}u^3+ 2321306389t^{12}u^4 -328391981t^{11}u^5\\
&+ 2250079387t^{10}u^6+ 1735063288t^9u^7-166127471t^8u^8+ 1710107151t^7u^9 -4967945152t^6u^{10}\\
&+ 2909685059t^5u^{11} -3065914682t^4u^{12} -738048656t^3u^{13}+ 1326185934t^2u^{14} -1625065825tu^{15}+ 1265968623u^{16}, \\
f_4 &= 4049930632t^{16}+ 2024965316t^{15}u -1469760896t^{14}u^2+ 3304162280t^{13}u^3 -2127737989t^{12}u^4+ 257670216t^{11}u^5\\
&+ 3908162477t^{10}u^6 -8926764895t^9u^7+ 4375575638t^8u^8 -2051118809t^7u^9+ 1008941202t^6u^{10}\\
&+ 6106860508t^5u^{11} -1555226010t^4u^{12}+ 657977093t^3u^{13} -906871421t^2u^{14} -1265968623tu^{15}.
\end{split}
}
\]
We now compute the revlex initial ideal of the image of the corresponding element of $M^{16}_4$, via the ``map" command, exactly as in \cite[Pf of Claim 5.3]{C2}. The output, namely
\[
{\small
\beg{split}
&(x_1^4, x_0x_1^3,x_0^2x_1^2,x_0^3x_1, x_0^4, x_0^3x_2^2, x_0^2x_1x_2^2, x_0x_1^2x_2^2, x_1^3x_2^2, x_1^3x_3^2, x_0x_1^2x_3^2, x_0^2x_1x_3^2, x_0^3x_3^2, x_2^4x_3, x_1x_2^3x_3,\\
&x_0x_2^3x_3, x_1^2x_2^2x_3, x_0x_1x_2^2x_3, x_0^2x_2^2x_3, x_1^3x_2x_3, x_0x_1^2x_2x_3, x_0^2x_1x_2x_3, x_0^3x_2x_3, x_2^5, x_1x_2^4, x_0x_2^4, x_1^2x_2^3, x_0x_1x_2^3, x_0^2x_2^3),
\end{split}
}
\]
has minimal generators in degree 5 and and below; it follows that $i=0$, as desired.

\subsection{Codimension estimates from singularities}
In this subsection, we complete the proof of the nonexistence of nondegenerate rational curves of degree 16 on a general heptic by proving the following technical result, along the same lines as \cite[Lem 3.4a]{JK} and \cite[Thm 5.1]{C2}. Here $M_{n,g}^d \sub M_n^d$ denotes the subspace of morphisms whose images are of arithmetic genus $g$.
\beg{lem}\label{technical_lemma}
$\mbox{cod}(M_{5,g}^{16},M_5^{16}) \geq \min(3g,9)$.
\end{lem}

\beg{proof}[Proof of Lemma~\ref{technical_lemma}]
We may assume $g=3$ without loss of generality. Rational curves of genus 3 are naturally stratified by partitions of 3, corresponding to (the delta-invariants of) their local singularity types. 
Now follow \cite[pp. 36-37]{C3}. 
The idea is to interpret spaces of equisingular deformations as degeneracy loci for natural evaluation maps to certain principal parts bundles constructed out of the restricted tangent bundle. Positivity of the latter ensures these degeneracy loci have the expected dimension.

{\fl \bf Step one: Enumeration of singularity types.} Given a singularity $\widehat{C}$ of embedding dimension $r$ with $l$ branches parameterized by
\[
t_i \mapsto (f_i^{(1)}(t_i), \dots, f_i^{(r)}(t_i)), 0 \leq i \leq l,
\]
there is a corresponding map
\[
\phi: \mb{C}[[x_1,\dots,x_r]] \ra \mb{C}[[t_1]] \times \cdots \times \mb{C}[[t_l]]
\]
of power series rings, given by
\[
1 \mapsto \underbrace{(1, \dots, 1)}_{l \text{ times}}, x_1 \mapsto (f_1^{(1)}(t_1), \dots, f_l^{(1)}(t_l)), \dots, x_r \mapsto (f_1^{(r)}(t_1), \dots, f_l^{(r)}(t_l)).
\]
The {\it delta-invariant} $\de=\de(\widehat{C})$ is equal to the colength of the subring generated by the image of $\phi$ inside the infinite-dimensional vector space underlying $\mb{C}[[t_1]] \times \cdots \times \mb{C}[[t_l]]$. A more geometric formulation, also useful, is that $\de$ computes the ``local number of nodes" to which $\wh{C}$ may be smoothed. 

{\fl We} now examine possible singularities with $\de=1,2,3$.

\beg{itemize}
\item $\de=1$: It is well-known that the only possibilities are a node (parameterization $t_1 \mapsto (t_1,0), t_2 \mapsto (0,t_2)$) or a cusp (parameterization $t \mapsto (t^2,t^3)$) 

\item $\de=2$: We sort possibilities according to the number of branches, as follows.
\beg{itemize}
\item {\bf One branch.} Up to isomorphism of analytic local rings, the only possibilities are associated to parameterizations $t \mapsto (t^3,t^4,t^5)$ and $t \mapsto (t^2,t^5)$.
\item {\bf Two branches.} Possibilities are parameterized by $t_1 \mapsto (t_1,0,0), t_2 \mapsto (0,t_2^2,t_2^3)$ (union of cusp and a line in $\mb{C}^3$) and $t_1 \mapsto (t_1,0), t_2 \mapsto (t_2,t_2^2)$ (tacnode). 
\item {\bf Three branches.} The unique possibility is the spatial triple point parameterized by $t_1 \mapsto (t_1,0,0), t_2 \mapsto (0,t_2,0), t_3 \mapsto (0,0,t_3)$.
\end{itemize}

\item $\de=3$. We again sort possibilities according to number of branches.
\beg{itemize}
\item {\bf One branch.} Possibilities are parameterized by $t \mapsto (t^4,t^5,t^6,t^7)$ and $t \mapsto (t^3,t^5)$.
\item {\bf Two branches.} Possibilities are parameterized by 
$t_1 \mapsto (0,t_1), t_2 \mapsto (t_2^2,t_2^3)$ (cusp and line in $\mb{C}^2$);  $t_1 \mapsto (t_1,0), t_2 \mapsto (t_2,t_2^3)$ (two smooth branches inflectionally tangent to one another); and $t_1 \mapsto (t_1^3,t_1^4,t_1^5,0), t_2 \mapsto (0,0,0,t_2)$ and $t_1 \mapsto (t_1^2,t_1^5,0), t_2 \mapsto (0,0,t_2)$ (each the union of a line and a unibranch singularity with $\de=2$).
\item {\bf Three branches.} Possibilities include a planar triple point, parameterized by $t_1 \mapsto (t_1,0), t_2 \mapsto (0,t_2), t_3 \mapsto (t_3,t_3)$; a spatial triple point, two of whose branches share a tangent line; and the union of a cusp and two lines in $\mb{C}^4$.
\item {\bf Four branches.} The unique possibility is the quadruple point in $\mb{C}^4$ parameterized by $t_1 \mapsto (t_1,0,0,0), t_2 \mapsto (0,t_2,0,0), t_3 \mapsto (0,0,t_3,0), t_4 \mapsto (0,0,0,t_4)$.
\end{itemize}
\end{itemize}

{\fl \bf Step two: Bounding deformations that fix singularity types.} Let $\psi: \mb{P}^1 \ra \mb{P}^n$ be a morphism with a single unibranch singularity; in local analytic coordinates near the singularity $p$, $\psi$ is of the form
\[
\psi: t \mapsto (\psi_1(t), \dots, \psi_n(t))
\]
where
\[
\psi_i(t)= \fr{t^{r_i}}{r_i!}+ (\text{higher order terms}), 2 \leq r_1 \leq \cdots \leq r_n.
\]
The {\it ramification type} of $\psi$ is the $n$-tuple $(r_1,\cdots,r_n)$.

{\fl Likewise}, an infinitesimal deformation of $\psi$ that fixes $p$ and $\psi(p)$ is locally of the form
\[
\psi^{\ep}(t)=(\psi_1(t)+ \ep(\al_1^{(1)}t+ \fr{\al_2^{(1)}}{2!}+ \cdots), \dots, \psi_n(t)+ \ep(\al_1^{(n)}t+ \fr{\al_2^{(n)}}{2!}+  \cdots))
\]
where $\al_j^{(i)}$ are complex numbers, for all $i$ and $j$, and $\ep^2=0$.

{\fl Those} deformations under which the ramification type of our singularity also stays fixed satisfy
\beg{equation}\label{diffcondition}
\text{rank}(\mbox{trunc}_{r_i-1}(\psi^{\ep})) \leq i-1, 1 \leq i \leq n
\end{equation}
where $\mbox{trunc}_j$ is the $j$th ``truncation" map that sends a power series in $t$ to the $j$-tuple of its first $j$ derivatives. 

{\fl The local condition \eqref{diffcondition} means exactly} that all $i \times i$ minors of the $n \times (r_i-1)$ matrix $A^{\psi}$
whose $(k,l)$th entry is
\[
A^{\psi}_{k,l}= \fr{d^l}{dt^l} \psi_k^{\ep}(t)
\]
vanish.


{\fl Note, moreover, that} every summand of $\psi^*T_{\mb{P}^n}$ has degree at least $\deg(\psi)+1$, as one sees by pulling back of the Euler sequence on $\mb{P}^n$, and using the fact that $\psi$ is nondegenerate. It follows that the canonical map
\[
\mbox{ev}: H^0(\mb{P}^1, \psi^*T_{\mb{P}^n}(-p)) \ra H^0(\mb{P}^1,\psi^* T_{\mb{P}^n}(-p))/\psi^* T_{\mb{P}^n}(-r_n p))
\]
is surjective whenever $r_n \leq \deg(\psi)+1$.

{\fl Consequently}, sections $\ka \in H^0(\mb{P}^1,\psi^*T_{\mb{P}^n})$ that define deformations preserving ramification type, $p$, and $\psi(p)$, i.e., whose images 
\[
\widetilde{\ka} \in H^0(\mb{P}^1,\psi^* T_{\mb{P}^n}(-p))/\psi^* T_{\mb{P}^n}(-r_n p))
\]
satisfy the determinantal conditions described in the preceding paragraph, comprise the pullback (via $\mbox{ev}$) of a Schubert cycle $\sig$ of codimension $\sum_{i=1}^n (r_i-i)$ in $G(n,H^0(\mb{P}^1,\psi^*T_{\mb{P}^n}(-p)/(-r_n p)))$. More precisely, letting
\[
\mc{I} \hra H^0(\mb{P}^1,\psi^*T_{\mb{P}^n}(-p)/(-r_n p)) \times G(n,H^0(\mb{P}^1,\psi^*T_{\mb{P}^n}(-p)/(-r_n p)))
\]
denote the incidence correspondence of sections contained in 5-dimensional subspaces, equipped with projections $\pi_1$ and $\pi_2$ onto its two factors, the deformations of interest comprise $\pi_1 (\pi_2^{-1} \sig)$. Since $\pi_1$ and $\pi_2$ are each equidimensional and surjective, it follows  that the locus of deformations is of the expected codimension $(\sum_{i=1}^n r_i -\binom{n}{2})$ in $H^0(\mb{P}^1,\psi^*T_{\mb{P}^n}(-p))$.

{\fl A similar} analysis applies to deformations that fix ramification types, and orders of contact, of singularities with multiple branches, such as those that appear in Step one. We give the details of this analysis in the case of the quadruple point in $\mb{C}^4$, and leave the remaining cases to the reader.

{\fl Namely}, consider the set of (infinitesimal) deformations that fix the preimages $p_1, p_2, p_3, p_4$ of the quadruple point, as well as their common image on the target. Such deformations correspond to sections of
\[
H^0(\mb{P}^1, \psi^* T_{\mb{P}^5}(-p_1-p_2-p_3-p_4)).
\]
Since $\psi^* T_{\mb{P}^5}$ is of rank 5 and each summand is of degree at least 17, it follows that such deformations define a codimension of at least 20 inside $M_4^{16}$. Varying the choice of $p_1, p_2, p_3, p_4$ and their common target $\psi(p_1)$, we obtain a locus of codimension at least 11, which beats the required estimate of 9.
\end{proof}

\subsection{Concluding remarks}
Unlike in \cite{C1} and \cite{C2}, I have not made a serious attempt at maximizing the degree $d$ for which the methods used above prove nonexistence of rational curves of degree $d$ on a general heptic, though I strongly suspect the methods can be pushed several degrees beyond 16. On the other hand, the heptic provides the first example of a hypersurface of {\it general type} for which a Clemens-type theorem is not known. It seems plausible (if not at all clear) that an adjunction-based method that exploits the positivity of the canonical bundle of the hypersurface along the lines developed in \cite{Cl}, \cite{P}, and \cite{Vo}, and more recently, in \cite{DMR} could be used to establish nonexistence in {\it every} degree. Our motivation in writing this note is rather to demonstrate that when the degree $d$ of the rational curves in question is sufficiently low, the method developed in \cite{C1} and \cite{C2} yields a relatively cheap proof of nonexistence. Furthermore, the paper \cite{HJ} only partially exploited the potential of the method, and it does not appear to treat the case of highly singular rational curves (i.e., those with large arithmetic genus). This paper should be considered, therefore, as a further advertisement of a method whose advantage is that it is elementary, and applies irrespective of the birational type of the hypersurfaces involved, but whose disadvantage is also precisely that it does not take this extra data into account. It should also be considered in some sense as an amplification and reworking of \cite{HJ}.

\beg{thebibliography}{25}
\bibitem{Ba} E. Ballico, {\it On singular curves in positive characteristic}, Math. Nachr. 141 (1989), 267-73.
\bibitem{C1} E. Cotterill, {\it Rational curves of degree 10 on a general quintic threefold}, Comm. Alg. {\bf 33} (2005), 1833--72.
\bibitem{C2} E. Cotterill, {\it Rational curves of degree 11 on a general quintic 3-fold}, Quart. J. Math., published electronically on 23 February 2011; doi:10.1093/qmath/har001.
\bibitem{C3} E. Cotterill, {\it Rational curves of degree 11 on a general quintic 3-fold}, \url{arXiv:0711.2758v1}.
\bibitem{Cl} H. Clemens, {\it Lower bounds on genera of subvarieties of generic hypersurfaces}, Comm. Alg. {\bf 31} (2003), 3673--3711.
\bibitem{DMR} S. Diverio, J. Merker, and E. Rousseau, {\it Effective algebraic degeneracy}, Invent. Math. {\bf 180} (2010), no. 1, 161--223.
\bibitem{GLP} L. Gruson, R. Lazarsfeld, and C. Peskine, {\it On a theorem of Castelnuovo and the equations defining space curves}, Invent. Math. {\bf 72} (1983), 491--506.
\bibitem{GP} L. Gruson and C. Peskine, {\it Courbes de l'espace projectif: vari\'{e}t\'{e}s de s\'{e}cantes}, in "Enumerative geometry and classical algebraic geometry (Nice, 1981)," Progr. Math. 24 (1982), 1-31.
\bibitem{Gr} M. Green, {\it Generic initial ideals}, in ``Six lectures on commutative algebra,'' J. Elias, J.M. Giral, R.M. Miro-Roig, and S. Zarzuela (Eds.), Birkh\"{a}user, Boston, 1998.
\bibitem{GS} D. Grayson and M. Stillman, {\it Macaulay2: a software system for research in algebraic geometry}, available at \url{http://www.math.uiuc.edu/Macaulay2}.
\bibitem{Ha} J. Harris, ``Curves in projective space," Les Presses de l'Universit\'{e} de Montr\'{e}al, Montr\'{e}al, Qu\'{e}bec, 1982.
\bibitem{HJ} G. Hana and T. Johnsen, {\it Rational curves on a general heptic fourfold}, Bull. Belg. Math. Soc. Simon Stevin {\bf 16} (2009), 861--85.
\bibitem{JK} T. Johnsen and S. Kleiman, {\it Rational curves of degree at most 9 on a general quintic threefold}, Comm. Alg. {\bf 24} (1996), 2721--53.
\bibitem{P} G. Pacienza, {\it Rational curves on general projective hypersurfaces}, J. Alg. Geom. {\bf 12} (2003), no. 2, 245--267.
\bibitem{R} L. Ramella, {\it La stratification du sch\'ema de Hilbert des courbes rationnelles par le fibr\'e tangent restreint}, C. R. Acad. Sci. Paris {\bf 311} (1990), 181--4.
\bibitem{V} J.L. Verdier, {\it Two-dimensional sigma-models and harmonic maps from $S^2$ to $S^n$}, in ``Group Theoretical Methods in Physics", Springer Lecture Notes in Physics {\bf 180} (1983), 136--41.
\bibitem{Vo} C. Voisin, {\it A correction to `A conjecture of Clemens on rational curves on hypersurfaces'}, J. Diff. Geom. {\bf 49} (1998), 601--611.
\end{thebibliography}
\end{document}